\documentclass[12pt,a4paper]{article}
\usepackage{amsfonts}
\usepackage{mathtext}
\usepackage[T2A]{fontenc}
\usepackage{amssymb}
\usepackage{fancyhdr}
\usepackage{cite}
\usepackage[cp866]{inputenc}
\usepackage[english]{babel}
\usepackage{bm}
\usepackage{amsmath}
\usepackage{hyperref} 
\usepackage[dvips]{graphicx}
\usepackage[a4paper,left=30mm,right=15mm,bottom=30mm,top=30mm]{geometry}
\makeatletter
\renewcommand\@biblabel[1]{#1.}
\makeatother

\begin{document}

\newcommand{\pp}{\leftthreetimes}
\renewcommand{\le}{\leqslant}
\renewcommand{\ge}{\geqslant}
\renewcommand{\O}{\varnothing}
\newcommand{\y}{\breve{у}}               

\large

\title{\bf Finite groups with $\mathbb P$-subnormal primary \\ cyclic subgroups}

\author{V.\,N. Kniahina and V.\,S. Monakhov}

\date{}

\maketitle

\begin{abstract}
A subgroup $H$ of a group $G$ is called $\mathbb P$-{\sl subnormal} in $G$ 
whenever either $H=G$ or there is a chain of subgroups 
$H=H_0\subset H_1\subset \ldots \subset H_n=G$ 
such that $|H_i:H_{i-1}|$~is a prime for all $i$.
In this paper, we study the groups in which all primary cyclic subgroups 
are $\mathbb P$\nobreakdash-\hspace{0pt}subnormal.

\end{abstract}


MSC2010 20D20, 20E34

\bigskip

{\small }

\bigskip

{\noindent\bf Introduction} 

\medskip

We consider finite groups only.
A.\,F. Vasilyev, T.\,I. Vasilyeva and V.\,N.~Tyu\-tyanov in \cite{VVTSMJ}
introduced the following definition. Let $\mathbb P$ be the set of all prime 
numbers. A subgroup $H$ of a group $G$ is called $\mathbb P$-{\sl subnormal} in $G$ 
whenever either $H=G$ or there is a chain  
$$
H=H_0\subset H_1\subset \ldots \subset H_n=G
$$
of subgroups such that $|H_i:H_{i-1}|$~is prime for all $i$. 
 
Let $|G|=p_1^{a_1}p_2^{a_2} \ldots p_k^{a_k}$, where
$p_1>p_2> \ldots >p_k$, $a_i\in \mathbb N$.
We say that $G$ has an ordered Sylow tower
of supersolvable type if there exist normal subgroups $G_i$ with
$$
1=G_0\subset G_1\subset G_2\subset \ldots \subset G_{k-1}\subset G_k=G,
$$ 
and where each factor $G_{i}/G_{i-1}$ is 
isomorphic to a Sylow $p_i$-subgroup of $G$ for all~$i$.
We denote by $\mathfrak D$ the class of all groups 
which have an ordered Sylow tower of supersolvable type.
It is well known that $\mathfrak D$ is a hereditary saturated formation.

In \cite{VVTSMJ} finite groups with $\mathbb P$-subnormal 
Sylow subgroups were studied. A group $G$ is called 
$\mbox {w}$-supersolvable  
if every Sylow subgroup of $G$
is $\mathbb P$-subnormal in $G$. Denote by $\mbox {w}\mathfrak U$ the class
of all $\mbox {w}$-supersolvable groups. Observe that  the class $\mathfrak U$ 
of all supersolvable groups is included into $\mbox {w}\mathfrak U$. 
In \cite{VVTSMJ}, the authors proved that the class $\mbox {w}\mathfrak U$
is a saturated hereditary formation; every  group in $\mbox {w}\mathfrak U$
possesses an ordered Sylow tower of supersolvable type;  all metanilpotent and 
all biprimary subgroups in $\mbox {w}\mathfrak U$ are supersolvable.

In \cite{VVTGGU} the following problem was proposed.

\medskip

\noindent {\sl To describe the groups in which all  
primary cyclic subgroups are $\mathbb P$\nobreakdash-\hspace{0pt}subnormal.}

\medskip

In this note we solve this problem. Denote by $\mathfrak X$ the class 
of groups whose primary cyclic subgroups are all 
$\mathbb P$\nobreakdash-\hspace{0pt}subnormal. It is easy to verify that  
$\mathfrak U\subset \mbox {w}\mathfrak U\subset \mathfrak X$.  

\medskip
 
{\bf Theorem.} 1.  {\sl A group $G\in \mbox {w}\mathfrak U$ if and only if  
$G$ possesses an ordered Sylow tower of supersolvable type and every  
biprimary subgroup of $G$ is supersolvable.} 

2. {\sl The class $\mathfrak X$ is a hereditary saturated formation.}

3. {\sl A group $G\in\mathfrak X $ if and only if $G$ possesses an 
ordered Sylow tower of supersolvable type and every  
biprimary subgroup of $G$ with cyclic Sylow subgroup is supersolvable.}

4. {\sl Every minimal non-$\mathfrak X$-group is a biprimary minimal 
non-supersolvable group whose non-normal Sylow subgroup is cyclic. } 

\section{Preliminary results}

We use the standart notation of \cite{Hup}. The set of all 
prime divisors
of $|G|$ is denoted $\pi (G)$. We write $[A]B$ for a semidirect product 
with a normal subgroup $A$. If $H$ is a subgroup of $G$, then 
$\mbox{Core}_GH=\bigcap _{x\in G}x^{-1}Hx$ is called the core of 
$H$~in~$G$. If a group $G$ contains a maximal subgroup $M$ with trivial
core, then $G$ is said to be primitive and $M$ is its primitivator.
We will use the following notation: 
$S_n$ and $A_n$ are the symmetric and the
alternating groups of degree $n$, $E_{p^t}$~is the elementary abelian group
of order $p^t$, $Z_m$~is the cyclic group of order $m$, $D_8$ is the 
dihedral group of order $8$,  
$Z(G)$, $\Phi (G)$, $F(G)$, $G^{\prime}$ is the center, the Frattini 
subgroup, the Fitting subgroup and the derived subgroup of $G$ 
respectively.

\medskip 

{\bf Lemma 1.}
{\sl Let $H$ be a subgroup of a solvable group $G$, and 
assume that $|G:H|$~is a prime number. Then $G/\emph {Core}_GH$ is supersolvable.}

\medskip

{\sc Proof.}
By hypothesis, $|G:H|=p$, where $p$ is a prime number. If $H=\mbox {Core}_GH$, 
then $G/H$~is
cyclic of order $p$ and $G/\mbox {Core}_GH$ is supersolvable, 
as required.  Assume now that $H\ne \mbox {Core}_GH$, i.\,e.,
$H$ is not normal in $G$. It follows that $G/\mbox {Core}_GH$ contains a 
maximal subgroup $H/\mbox {Core}_GH$ with trivial core. 
Hence, $G/\mbox {Core}_GH$ is primitive and the Fitting subgroup 
$F/\mbox {Core}_GH$ of $G/\mbox {Core}_GH$ has prime order $p$. 
Since 
$$
F/\mbox {Core}_GH=C_{G/\mbox {Core}_GH}(F/\mbox {Core}_GH),
$$
it follows that 
$$
(G/\mbox {Core}_GH)/(F/\mbox {Core}_GH)\simeq H/\mbox {Core}_GH
$$ 
is isomorphic to a cyclic group of order dividing $p-1$. Thus, 
$G/\mbox {Core}_GH$ is supersolvable. 

\medskip

{\bf Lemma 2.}  1. {\sl A group $G$ is supersolvable if and only if
all of its maximal subgroups have prime indices.

$2$.  Every subgroup of a supersolvable group is $\mathbb P$-subnormal. }

\medskip

{\sc Proof.} 1. This is Huppert's classic 
theorem, see \cite[Theorem VI.9.5] {Hup}.

2. The statement follows from 1 of the lemma.

\medskip

Immediately, using the definition of $\mathbb P$-subnormality, 
we deduce the following properties.

\medskip

{\bf Lemma 3.} 
{\sl 
Suppose that $H$~is a subgroup of $G$, and let 
$N$~be a normal subgroup of $G$. Then the following hold:

$1)$ if $H$ is $\mathbb P$-subnormal in $G$, then 
$(H\cap N)$ is $\mathbb P$-subnormal in $N$, and $HN/N$ is
$\mathbb P$-subnormal in $G/N$;

$2)$ if $N\subseteq H$ and $H/N$ is $\mathbb P$-subnormal in $G/N$, 
then $H$ is $\mathbb P$-subnormal in~$G$;

$3)$ if $H$ is $\mathbb P$-subnormal in $K$, and $K$ is $\mathbb P$-subnormal in 
$G$, then $H$ is $\mathbb P$\nobreakdash-\hspace{0pt}subnormal in $G$;

$4)$ if $H$ is $\mathbb P$-subnormal in $G$, then $H^g$ is $\mathbb P$-subnormal in 
$G$ for each element $g\in G$}.

\medskip

{\bf Example 1.} The subgroup $H=A_4$ of the alternating group 
$G=A_5$  is $\mathbb P$-subnormal.
If $x\in G\setminus H$, then $H^x$ is $\mathbb P$-subnormal in $G$.
The subgroup $D=H\cap H^x$ is a Sylow 3-subgroup of $G$ and $D$
is not $\mathbb P$\nobreakdash-\hspace{0pt}subnormal in $H$. 
Therefore, an intersection of two $\mathbb P$-subnormal subgroups is not 
$\mathbb P$-subnormal.
Moreover, if $H$ is $\mathbb P$-subnormal in $G$ and 
$K$~is an arbitrary subgroup of $G$, in general, their intersection 
$H\cap K$ is not $\mathbb P$-subnormal in~$K$. 

\medskip

However, this situation is impossible if $G$ is a solvable group.

\medskip

{\bf Lemma 4.} 
{\sl 
Let $G$ be a solvable group. Then the following hold:

$1)$ if $H$ is $\mathbb P$-subnormal in $G$, and $K$ is a subgroup of $G$, 
then $(H\cap K)$ is $\mathbb P$-subnormal in $K$;

$2)$ if $H_i$ is $\mathbb P$-subnormal in $G$, $i=1,2$, then $(H_1\cap H_2)$ 
is $\mathbb P$-subnormal in~$G$.}

\medskip

{\sc Proof.} 1. It is clear that in the case $H=G$ the statement is true. Let 
$H\ne G$. According to the definition of $\mathbb P$-subnormality, there exists 
a chain of subgroups
$$
H=H_0\subset H_1\subset \ldots \subset H_n=G
$$
such that $|H_i:H_{i-1}|$~is a prime number for any $i$. 
We will use induction~by~$n$.  

Consider the case when $n=1$. In this situation, $H=H_{n-1}$~
is a maximal subgroup of prime index in $G$. By Lemma 1,
$G/N$ is supersolvable, $N=\mbox {Core}_GH$. Since, by Lemma 2~(2),
every subgroup of a supersolvable group is 
$\mathbb P$\nobreakdash-\hspace{0pt}subnormal, we have 
$$
H/N\cap KN/N=N(H\cap K)/N
$$
is $\mathbb P$-subnormal in $KN/N$. Lemma 3~(2) implies that 
$N(H\cap K)$ is $\mathbb P$-subnormal in $KN$. It means that there exists 
a chain of subgroups
$$
N(H\cap K)=A_0\subset A_1\subset \ldots \subset A_{m-1}\subset A_m=NK
$$
such that $|A_i:A_{i-1}|\in \mathbb P$ for all $i$. Since
$$
N(H\cap K)\subseteq A_i\subseteq  NK,
$$
we have $A_i=N(A_i\cap K)$  and
$H\cap K\subseteq  A_i\cap K$ for all $i$. We introduce the 
notation $B_i=A_i\cap K$. It is clear that 
$$
B_{i-1}\subseteq B_i, \ A_i=N(A_i\cap K)=NB_i, \ 
N\cap B_i=N\cap A_i\cap K=N\cap K 
$$
for all $i$. Since $N\subseteq H$, we have 
$$
B_0=A_0\cap K=N(H\cap K)\cap K=(N\cap K)(H\cap K)=H\cap K,  
$$
$$
B_m=A_m\cap K=KN\cap K=K.
$$
Moreover,
$$
|A_i:A_{i-1}|=\frac {|NB_i|}{|NB_{i-1}|}=\frac 
{|N||B_i||N\cap B_{i-1}|}
{|N\cap B_i||N||B_{i-1}|}=\frac{|B_i:B_{i-1}|}{|N\cap B_i:N\cap B_{i-1}|}=
$$
$$
=\frac{|B_i:B_{i-1}|}{|N\cap K:N\cap K|}= |B_i:B_{i-1}|.
$$
Now we have a chain of subgroups
$$
H\cap K=B_0\subset B_1\subset \ldots \subset B_{m-1}\subset 
B_m=K, \ |B_i:B_{i-1}|\in \mathbb P, \ 1\le i\le m,
$$
which proves that the subgroup $H\cap K$ is $\mathbb P$-subnormal in $K$. 

Let $n>1$. Since $H_{n-1}$~is a maximal subgroup of prime index in 
$G$, and $G$ is solvable, thus, as it was proved, $H_{n-1}\cap K$ is   
$\mathbb P$-subnormal in $K$. The subgroup $H$ is 
$\mathbb P$-subnormal in the solvable group $H_{n-1}$ and the induction is
applicable to them. By induction, 
$$
H\cap (H_{n-1}\cap K)=H\cap K
$$ 
is $\mathbb P$-subnormal in $H_{n-1}\cap K$. By Lemma 3~(3),  
$H\cap K$ is $\mathbb P$-subnormal in $K$. 
                               
2. Let $H_i$ is $\mathbb P$-subnormal in $G$, $i=1,2$. It follows from 1 of the 
lemma, that $(H_1\cap H_2)$ is $\mathbb P$-subnormal in $H_2$. Now by Lemma 3~(3), 
we obtain that  $(H_1\cap H_2)$ is $\mathbb P$-subnormal in $G$. 

\medskip 

{\bf Lemma 5.}
{\sl Let $H$ be a subnormal subgroup of a solvable group $G$. Then 
$H$ is $\mathbb P$-subnormal in $G$.}

\medskip

{\sc Proof.} Since $H$ is subnormal in $G$, and $G$ is solvable,   
then there exists a series 
$$
H=H_0\subset H_1\subset \ldots \subset H_{n-1}\subset H_n=G,
$$
such that $H_i$ is normal in $H_{i+1}$ for all $i$.
Working by induction on $|G|$, we can assume that $H$ is $\mathbb P$-subnormal 
in $H_{n-1}$.
Since $G/H_{n-1}$ is solvable, the composition factors of $G/H_{n-1}$
have prime orders. Thus, there is a chain of subgroups
$$
H_{n-1}=G_0\subset G_1\subset \ldots \subset G_{m-1}\subset G_m=G
$$
such that $G_j$ is normal in $G_{j+1}$ and $|G_{j+1}/G_j|\in \mathbb P$ for all 
$j$. This means that $H_{n-1}=G_0$ is $\mathbb P$-subnormal in $G$.
Using Lemma 3 (3), we deduce that $H$ is a $\mathbb P$-subnormal
in $G$.

\medskip  

{\bf Example 2.} The subgroup $Z(SL(2,13))$ 
of the non-solvable group $SL(2,13)$  
is normal, but is not $\mathbb P$-subnormal. This follows from the fact that
the identity subgroup is not 
$\mathbb P$-subnormal in $PSL(2,13)=SL(2,13)/Z(SL(2,13))$.

\medskip

{\bf Lemma 6.} {\sl Let $A$ be a $p$-subgroup of a group $G$. Then $A$ is 
subnormal in $G$ if and only if $A\subseteq O_p(G)$.}

\medskip

{\sc Proof.} The statement follows from Theorem 2.2\cite{Is}.

\medskip

{\bf Lemma 7.}
{\sl Let  $A$ be a $p$-subgroup of a group $G$. 
If $|G:N_G(A)|=p^\alpha $, $\alpha \in \mathbb N$, then $A$ is subnormal in $G$.}
                                                           
\medskip

{\sc Proof.}  Let $P$ be a Sylow $p$-subgroup of $G$ with the property that
$P$ contains $A$. Then
$$
G=N_G(A)P, \ A^G=A^{N_G(A)P}=A^{P}\subseteq P,
$$
so $A^G\subseteq O_p(G)$. It is clear that $A$ is subnormal in $G$. 

\medskip

{\bf Lemma 8.}
{\sl Let $p$ be the largest prime divisor of $|G|$, and let $A$~be a 
$p$\nobreakdash-\hspace{0pt}subgroup of $G$. If $A$ is $\mathbb P$-subnormal in $G$, then 
$A$ is subnormal in $G$.}
                                                           
\medskip

{\sc Proof.}  Let $|A|=p^\alpha$. Since $A$ is 
$\mathbb P$-subnormal in $G$, then there exists a series
$$
A=A_0\subset A_1\subset \ldots \subset A_{t-1}\subset A_t=G, \
|A_i:A_{i-1}|\in \mathbb P, \ 1\le i\le t.
$$
Since $|A_1:A_0|\in \mathbb P$, we have 
$$
|A_1|=p^{1+\alpha} \ \ or \ \ |A_1|=p^\alpha q, \ \ p\ne q.
$$
If $|A_1|=p^{1+\alpha}$, then $A$ is a normal subgroup of $A_1$. If 
$|A_1|=p^\alpha q$, then $p>q$ and again $A$ is normal in $A_1$. 
Suppose we already know that $A$ is subnormal in $A_j$. 
Using Lemma 6 we have, $A\subseteq O_p(A_j)$. Since $|A_{j+1}:A_j|\in \mathbb P$, 
we obtain 
$$
|A_{j+1}|=p|A_j| \ \  or \ \ |A_{j+1}|=q|A_j|, \ \  p\ne q.
$$
If $|A_{j+1}|=p|A_j|$, then, by Lemma 7, $O_p(A_j)\subseteq O_p(A_{j+1})$, and $A$
is subnormal in  $A_{j+1}$. If $|A_{j+1}|=q|A_j|$, $p\ne q$, then $p>q$. 
Consider the set of left cosets of $A_j$ in $A_{j+1}$. 
We know that $A_{j+1}/\mbox{Core}_{A_{j+1}} A_j$ is isomorphic to a subgroup of the 
symmetric group $S_q$ and any Sylow $p$-subgroup of $A_{j+1}$ is contained in 
$\mbox{Core}_{A_{j+1}}A_j$. Since $A$ is subnormal in $A_j$,
so $A$ is subnormal in $\mbox{Core}_{A_{j+1}} A_j$. Since 
$\mbox{Core}_{A_{j+1}} A_j$ is normal in 
$A_{j+1}$, it follows that $A$ is subnormal in  $A_{j+1}$. 
Therefore, $A$ 
is subnormal in $A_i$ for all $i$. This implies that $A$ is subnormal~in~$G$. 

\medskip

{\bf Corollary.} (\cite[Proposition 2.8]{VVTSMJ}) {\sl 
Every $\mbox w$-supersolvable group possesses an ordered Sylow tower 
of supersolvable type.}

\medskip

{\sc Proof.}  Use induction on $|G|$. Let $G$~be a $\mbox w$-supersolvable 
group, and assume that $p$~is the largest prime divisor of $|G|$.
Let $P$ be a Sylow $p$-subgroup of $G$. By Lemma 8, $P$ is normal 
in $G$.  It follows by Lemma 3~(1), that any quotient
group of a $\mbox w$-supersolvable group is $\mbox w$-supersolvable. 
Working by induction on $|G|$, we deduce that $G/P$ possesses an ordered Sylow tower 
of supersolvable type, so $G$ possesses an ordered Sylow tower 
of supersol\-vable type.  The corollary is proved.

\medskip 

Recall that a Schmidt group is a finite non-nilpotent group all of whose 
proper subgroups are nilpotent.  Given a class $\mathfrak F$~of 
groups. By $\mathcal M (\mathfrak F)$ we denote the class of all minimal 
non\nobreakdash-\hspace{0pt}$\mathfrak F$\nobreakdash-\hspace{0pt}groups. 
A group $G$ is a minimal 
non\nobreakdash-\hspace{0pt}$\mathfrak F$\nobreakdash-\hspace{0pt}group if 
$G\notin\mathfrak F$ but all proper subgroups of $G$
belong to $\mathfrak F$. Clearly, the class $\mathcal M (\mathfrak N)$ 
consists of Schmidt groups. Here $\mathfrak N$~denotes the class of all nilpotent
groups. We will need the properties of groups from  
$\mathcal M (\mathfrak N)$ and $\mathcal M (\mathfrak U)$.

\medskip

{\bf Lemma 9.} {(\cite{MonUMK}, \cite{Red})} 
 {\sl 
Let $S\in \mathcal M (\mathfrak N)$. Then the following statements hold:

$1)$ $S=[P]\langle y \rangle$, where $P$ is a normal Sylow $p$-subgroup,
and $\langle y \rangle$~ is a non-normal cyclic Sylow $q$-subgroup,
$p$ and $q$ are distinct primes, ${y^q\in Z(S)}$;

$2)$ $|P/P^{\prime}|=p^m$, where $m$ is the order of $p$ modulo $q$;
                                                      
$3)$ if $P$ is abelian, then $P$  is an elementary abelian $p$-group 
of order $p^m$ and $P$ is a minimal normal subgroup of $S$;  

$4)$ if $P$ is non-abelian, then $Z(P)=P^{\prime}=\Phi (P)$ and
${|P/Z(P)| = p^m}$;

$5)$ $Z(S)=\Phi (S)=\Phi (P)\times \langle y^q\rangle$;
$S^{\prime}=P$, $P^{\prime}=(S^{\prime})^{\prime}=\Phi (P)$;

$6)$ if $N$ is a proper normal subgroup of $S$, then
$N$ does not contain $\langle y \rangle$ and 
either $P\subseteq N$ or $N\subseteq \Phi (S)$.}

\medskip

{\bf Lemma 10.} (\cite{Doerk})
{\sl Let $G\in \mathcal M (\mathfrak U)$. 
Then the following statements hold:

$1)$ $G$ is solvable and $|\pi (G)|\le 3$;

$2)$ if $G$ is not a Schmidt group, then $G$ possesses an ordered 
Sylow tower of supersolvable type;

$3)$ $G$ has a unique normal Sylow subgroup $P$
and $P=G^{\frak U}$;

$4)$ $|P/\Phi(P)|>p$ and $P/\Phi (P)$ is a minimal normal subgroup    
of~$G/\Phi (G)$;

$5)$ the Frattini subgroup $\Phi (P)$ of $P$ is supersolvable 
embedded in $G$, i.e., there exists a series
$$
1\subset N_0\subset N_1 \ldots \subset  N_m=\Phi (P)
$$
such that $N_i$ is a normal subgroup of $G$ and $|N_i/N_{i-1}|\in \mathbb P$ for
all $i$;

$6)$ let $Q$ be a complement to $P$ in $G$, then $Q/Q\cap \Phi (G)$~is
a minimal non-abelian group or a cyclic group of prime power order;

$7)$ all maximal subgroups of non-prime index are conjugate in $G$,
and moreover, they are conjugate to $\Phi (P)Q$.}
                      
\medskip   

We now present new properties of $\mbox w$-supersolvable groups.

\medskip    

{\bf Lemma 11.}
1. {\sl If $G\in \mathcal M(\mathfrak U)$ and $|\pi (G)|=3$, then $G$ is
$\emph w$-supersolvable.}

2. $\mathcal M(\mathfrak U)\setminus \emph w\mathfrak U=\{G\in \mathcal 
M(\mathfrak U)\mid |\pi (G)|=2\}$.
                                                           
3. {\sl If $G\in \emph w\mathfrak U$,  
then  the derived length 
of $G/\Phi (G)$ is at most $|\pi (G)|$.}

\medskip 

{\sc Proof.} 1. Let $G\in \mathcal M(\mathfrak U)$ and $|\pi (G)|=3$.
By Lemma 10, $G=[P]([Q]R)$, where $P$, $Q$ and $R$ are Sylow 
subgroups of $G$. The subgroup $P$ is normal in $G$, and using Lemma 5,  
we see that $P$ is $\mathbb P$-subnormal in $G$. The subgroup $PQ$ is normal
in $G$, and by Lemma 5, $PQ$ is $\mathbb P$-subnormal in $G$. Since $PQ$
is supersolvable, it follows by Lemma 2~(2), that $Q$ 
is $\mathbb P$-subnormal in $PQ$. By Lemma 3~(3), $Q$ is 
$\mathbb P$-subnormal in  $G$. Since $G/P\simeq QR$, so $G/P$ 
is supersolvable and $PR/P$ is $\mathbb P$-subnormal in $G/P$ by Lemma 2~(2). 
Hence $PR$ is $\mathbb P$-subnormal in $G$ by Lemma 3~(2).
Since $PR$ is supersolvable, we see that $R$ is $\mathbb P$-subnormal in $PR$
by Lemma 2~(2). Hence $R$ is $\mathbb P$-subnormal in $G$ by Lemma 3~(3).
We conclude that all Sylow subgroups of $G$ are $\mathbb P$-subnormal in $G$. 
Therefore, $G$ is $\mbox w$-supersolvable.

2. If $G\in \mathcal M(\mathfrak U)\setminus \mbox w\mathfrak 
U$, then $|\pi (G)|=2$ by assertion 1 of the lemma. Conversely, let 
$G\in \{\mathcal M(\mathfrak U)\mid |\pi (G)|=2\}$. Suppose that 
$G\in \mbox w\mathfrak U$. Then by Theorem 2.13~(2) \cite{VVTSMJ}, 
the group $G$ is supersolvable, this is a contradiction. 

3. By theorem 2.13~(3) \cite{VVTSMJ}, 
$G/F(G)$ has only abelian Sylow subgroups.
By theorem VI.14.16 \cite{Hup}, the derived length
of $G/F(G)$ is at most  $|\pi (G/F(G))|$. Since $G$ has  
an ordered Sylow tower of supersolvable type, so
$|\pi (G/F(G))|\le |\pi (G)|-1$. But if $G$ is solvable, then
the quotient group $F(G)/\Phi (G)$ is abelian, and we conclude that 
the derived length of $G/\Phi (G)$ does not exceed $|\pi (G)|$. 

\bigskip

\section{Finite groups with $\mathbb P$-subnormal primary cyclic subgroups}                       

\medskip

{\bf Example 3.}
There are three non-isomorphic minimal non-supersolvable groups of order 400: 
$$
[E_{5^2}](<a><b>), \ |a|=|b|=4.
$$
Numbers of these groups in the library of SmallGroups~\cite{GAP}
are [400,129], [400,130], [400,134]. Sylow 2-subgroups of these groups are 
non-abelian and have the form: $[Z_4\times Z_2]Z_2$ and $[Z_4]Z_4$. 
Let $G$ be one of these groups. All subgroups
of $G$ are $\mathbb P$-subnormal except the maximal subgroup 
$<a><b>$.  Therefore, these groups belong to the class $\mathfrak X$.

\medskip 

{\bf Example 4.}
The general linear group $GL(2,7)$ contains the symmetric group $S_3$
which acts irreducibly on the elementary abelian group $E_{7^2}$ of order 49. 
The semidirect product $[E_{7^2}]S_3$ is a minimal non-supersolvable 
group, it has subgroups of orders 14 and 21.
Every primary cyclic subgroup of the group $[E_{7^2}]S_3$ is $\mathbb P$-subnormal. 
Therefore, these group belong to the class $\mathfrak X$.

\medskip 

{\bf Example 5.} Non-supersolvable Schmidt groups do not belong 
to the class $\mathfrak X$. We verify this fact. 
Let $S=[P]Q$~be a non-supersolvable Schmidt group. 
Suppose that $S\in \mathfrak X$. It follows that $Q$ is $\mathbb P$-subnormal 
in $S$ and $Q$ is contained in some subgroup $M$ of prime index. 
Therefore, $M=P_1\times Q$, where $P_1$~is a normal subgroup of $S$ 
with the property 
$|P/P_1|=p$. By the properties of Schmidt groups, see Lemma 9,
we have $|P/\Phi (P)|>p$ and $P/\Phi (P)$~ is a chief factor of $S$.
We have a contradiction.                            

\medskip

{\bf Lemma 12.}
{\sl Suppose that all cyclic $p$-subgroups of a group $G$ are $\mathbb 
P$\nobreakdash-\hspace{0pt}subnormal and let $N$~be a normal subgroup of $G$. 
Then all cyclic subgroups of $N$ and $G/N$ are $\mathbb P$-subnormal.}
                                                           
\medskip 

{\sc Proof.}  Lemma 3 (1) implies that all cyclic $p$-subgroups of the 
normal group $N$ are $\mathbb P$-subnormal in $N$. Let $A/N$~be a cyclic
$p$-subgroup of $G/N$ and assume that $a\in A\setminus N$. 
Let $P$~be a Sylow $p$-subgroup of $\langle a\rangle$. By hypothesis, 
$P$ is $\mathbb P$-subnormal in $G$. Since $PN/N=AN/N$, it follows by Lemma 3~(1), 
that $A/N$ is $\mathbb P$-subnormal in $G/N$.

\medskip

{\bf Lemma 13.}
1. {\sl If every primary cyclic subgroup of a group $G$ is 
$\mathbb P$\nobreakdash-\hspace{0pt}subnormal, then $G$ possesses an ordered 
Sylow tower of supersolvable type.}

2. $\mathfrak U\subset \mbox {w}\mathfrak U\subset \mathfrak  X\subset \mathfrak D$.
                                                          
\medskip

{\sc Proof.}  1. Let $P$~be a Sylow $p$-subgroup of $G$, where $p$ is the largest
prime divisor of $|G|$. If $a\in P$, then by hypothesis, the subgroup
$\langle a\rangle$ is $\mathbb P$\nobreakdash-\hspace{0pt}subnormal in $G$. 
By Lemma 8, the subgroup $\langle a\rangle$ is subnormal in $G$, and by Lemma 6, 
$\langle a\rangle\subseteq O_p(G)$.
Since $a$~is an arbitrary element of $P$, we see that $P\subseteq O_p(G)$, 
and hence $G$ is $p$-closed. By Lemma 12, the conditions of the lemma are  
inherited by all quotient groups of $G$. Applying induction on $|G|$,
we see that $G/P$ possesses an ordered Sylow tower of supersolvable type, and thus 
$G$ possesses an ordered Sylow tower of supersolvable type. 

2. By Lemma 2~(2), we have the inclusion $\mathfrak U\subseteq 
\mbox {w}\mathfrak U$. It follows from Example 4, that  
$[E_{7^2}]S_3$ is non-supersolvable and $[E_{7^2}]S_3\in \mbox {w}\mathfrak U
\setminus \mathfrak U$. Therefore, $\mathfrak U\subset \mbox {w}\mathfrak U$. 

We verify the inclusion $\mbox {w}\mathfrak U \subseteq \mathfrak X$. 
Suppose that $G\in \mbox {w}\mathfrak U$, and let
$A$~be an arbitrary primary cyclic subgroup of $G$. Then $A$ 
is a $p$-subgroup for some $p\in \pi (G)$. By Sylow's theorem,  
$A$ is contained in some Sylow $p$\nobreakdash-\hspace{0pt}subgroup $P$ 
of the group $G$. 
Since $G\in \mbox {w}\mathfrak U$, it follows that $P$ is $\mathbb P$-subnormal
in $G$. By Lemma 2~(2), $A$ is $\mathbb P$-subnormal in $P$,
and by Lemma 3~(3), $A$ is $\mathbb P$-subnormal in $G$. 
Therefore, $G\in \mathfrak X$. The group $[E_{5^2}]Q$ from Example 3 is a
biprimary minimal non-supersolvable group, $Q$ is non-cyclic.
The group $[E_{5^2}]Q\in \mathfrak X\setminus \mbox {w}\mathfrak U$,
therefore, $\mbox {w}\mathfrak U\subset \mathfrak X$. 

By the above assertion of the lemma,  
$\mathfrak  X\subseteq \mathfrak D$. 
Since there exist non-supersolvable Schmidt groups which 
have an ordered Sylow tower of supersolvable type
(for example, $[E_{5^2}]Z_3$), and they do not belong to the class 
$\mathfrak X$, it follows that $[E_{5^2}]Z_3\in \mathfrak D\setminus 
\mathfrak X$. Therefore, $\mathfrak  X\subset \mathfrak D$.

\medskip 

{\bf Lemma 14.}
{\sl Let $G$ be a minimal non-supersolvable group.  The group 
${G\not \in \mathfrak X}$ if and only if $G$ is a biprimary group whose  
non-normal Sylow subgroup is cyclic.}
                                                           
\medskip 

{\sc Proof.} 
Let $G\in \mathcal M(\mathfrak U)\setminus \mathfrak X$. 
If $|\pi (G)|=3$, then by Lemma 11~(1), $G\in \mbox w\mathfrak U$. Since
$\mbox w\mathfrak U\subset \mathfrak X$, we have $G\in \mathfrak X$, 
which contradicts the choice of $G$. So, if 
$G\in \mathcal M(\mathfrak U)\setminus \mathfrak X$, then $|\pi (G)|=2$ 
and $G=[P]Q$, where $P$ is a Sylow $p$-subgroup of $G$, $Q$ is a Sylow 
$q$-subgroup of $G$. Suppose that $Q$ is non-cyclic, and let $a\in Q$. Since 
$P\langle a\rangle$~is a proper subgroup of $G$, we deduce that 
$P\langle a\rangle$ is supersolvable. 
Lemma 2~(2) implies that $\langle a\rangle$ is a $\mathbb 
P$-subnormal subgroup of $P\langle a\rangle$. Since $P\langle a\rangle$ 
is subnormal in $G$, it follows by Lemma 5, that $P\langle a\rangle$ is a
$\mathbb P$-subnormal subgroup of $G$. Now by Lemma 3~(3), we deduce that
$\langle a\rangle$ is  $\mathbb P$-subnormal in $G$. Applying 
Lemma 3~(4), we can conclude that all cyclic $q$-subgroups of $G$ are
$\mathbb P$-subnormal in $G$. Lemma 5 implies that all cyclic $p$-subgroups of 
$G$ are $\mathbb P$-subnormal in $G$. Thus $G\in \mathfrak X$. 
We have a contradiction. Therefore, the assumption is false and $Q$ is 
cyclic. 

Conversely, let $G\in \mathcal M(\mathfrak U)$, $|\pi (G)|=2$
and a non-normal Sylow subgroup $Q$ of $G$ is cyclic. Assume that 
$G\in \mathfrak X$. This implies that $Q$ is $\mathbb P$-subnormal in $G$, and so 
both Sylow subgroups of the group $G$ are $\mathbb P$-subnormal.  
Now, by Theorem 2.13~(2) \cite{VVTSMJ}, $G$ is supersolvable,
which is a contradiction. 

\medskip

{\bf Lemma 15.} (\cite{Baer}) {\sl 
If $P$ is a normal Sylow subgroup of a group $G$, 
then $\Phi (P)=\Phi(G)\cap P$.}

\medskip
 
{\bf Theorem.} 1.  {\sl A group $G\in \mbox {w}\mathfrak U$ if and only if  
$G$ possesses an ordered Sylow tower of supersolvable type and every  
biprimary subgroup of $G$ is supersolvable.} 

2. {\sl The class $\mathfrak X$ is a hereditary saturated formation.}

3. {\sl A group $G\in\mathfrak X $ if and only if $G$ possesses an 
ordered Sylow tower of supersolvable type and every  
biprimary subgroup of $G$ with cyclic Sylow subgroup is supersolvable.}

4. {\sl Every minimal non-$\mathfrak X$-group is a biprimary minimal 
non-supersolvable group whose non-normal Sylow subgroup is cyclic. } 

\medskip
 
{\sc Proof.} 1. If a group $G\in \mbox {w}\mathfrak U$, then $G$ 
possesses an ordered Sylow tower of supersolvable type by the corollary of 
Lemma 8, and every biprimary subgroup of $G$ is 
$\mbox w$-supersolvable by Lemma 4\,(1). We conclude 
by Lemma 10\,(2), that every  
biprimary subgroup of $G$ is supersolvable.

Conversely, suppose that a group $G$ possesses an ordered Sylow tower of 
supersolvable type and every biprimary subgroup of $G$ is 
supersolvable. Assume that $G$ is not $\mbox w$-supersolvable. Let us choose 
among all such groups a group with the smallest possible
$|\pi (G)|$. Then $|\pi (G)|\ge 3$ and $G$ contains a Sylow 
$r$-subgroup $R$ such that $R$ is not $\mathbb P$-subnormal in $G$. 
Let $p\in \pi (G)$, where $p$~is the largest prime divisor of $|G|$
and let $P$~be a Sylow $p$-subgroup of $G$. Since
$G$ possesses an ordered Sylow tower of supersolvable type, we deduce that 
$P$ is normal in $G$. By hypothesis, 
the subgroup $PR$ is supersolvable, and we deduce by Lemma 2\,(2) that 
$R$ is $\mathbb P$\nobreakdash-\hspace{0pt}subnormal in $PR$. It is clear that $G/P$ 
possesses an ordered Sylow tower of supersolvable type and all of its 
biprimary subgroups are supersolvable. Since $|\pi (G/P)|=|\pi (G)|-1$,
it follows by the inductive hypothesis, that $G/P$ is $\mbox w$-supersolvable. 
Therefore, the Sylow $r$-subgroup
$PR/P$ is $\mathbb P$-subnormal in $G/P$. Lemma 3\,(2) implies that  
the subgroup $PR$ is $\mathbb P$-subnormal in $G$, and hence by 
Lemma 3\,(3), the subgroup $R$ is $\mathbb P$-subnormal in $G$. 
This is a contradiction. 

2. By Lemma 13~(1), the class $\mathfrak X$ consists of finite groups
which have an ordered Sylow tower of supersolvable type, so we can apply Lemma 4. 
Let $G\in \mathfrak X$ and suppose that $H$~is an arbitrary subgroup of $G$. 
If $A$ is a cyclic primary subgroup of $H$, then $A$ is $\mathbb P$-subnormal 
in $G$. By Lemma 4~(1), the subgroup $A$ is $\mathbb P$-subnormal in $H$,
and hence $\mathfrak X$ is a hereditary class.

By Lemma 12, the class $\mathfrak X$ is closed under homomorphic image.
By induction on the order of $G$, we verify that the class $\mathfrak X$ 
is closed under subdirect products.  Let $G$~be a group 
of least order with the following properties:  
$$
G/N_i\in \mathfrak X, \ i=1,2, \ N_1\cap N_2=1, \ G\not \in \mathfrak X.
$$
In this case, $G$ has a primary cyclic subgroup $A$ 
which is not $\mathbb P$-subnormal in $G$. Since $G/N_i\in \mathfrak X$, $i=1,2$, it 
follows that $AN_i/N_i$ is $\mathbb P$-subnormal in $G/N_i$, and thus by 
Lemma 4~(2), $AN_1\cap AN_2$ is $\mathbb P$-subnormal in $G$. 
If $K=AN_1N_2$~is a proper subgroup of $G$, then $K/N_i\in \mathfrak X$ 
because $\mathfrak X$~is a hereditary class. By the induction hypothesis, 
$K\in \mathfrak X$. It follows that $A$ is $\mathbb P$-subnormal in 
$AN_1\cap AN_2$, and by Lemma 3~(3), $A$ is $\mathbb P$-subnormal in $G$, 
which is a contradiction. Therefore, $G=AN_1N_2$. Assume that $G=AN_1$. Then 
$$
N_2\simeq N_1N_2/N_1\subseteq   G/N_1\simeq A/A\cap N_1.
$$
Thus $N_2$ is cyclic and $AN_2$ is supersolvable by Theorem VI.10.1 \cite{Hup}.
It follows that $A$ is $\mathbb P$-subnormal in $AN_2$, $AN_2$ is 
$\mathbb P$-subnormral in $G$, and by Lemma 3~(3), $A$ 
is $\mathbb P$-subnormal in $G$, which is a contradiction.
Thus our assumption is false and $AN_1\ne G\ne AN_2$.

The subgroup $D=N_1\cap AN_2$ is normal in $AN_2=H\ne G$.
Hence the group $H$ contains two normal subgroups $D$ and $N_2$ such that    
$$
D\cap N_2\subseteq N_1\cap N_2=1,
$$
$$
H/N_2\subset G/N_2\in \mathfrak X, \ H/N_2\in \mathfrak X,
$$
$$
G/N_1=(AN_2)N_1/N_1\simeq AN_2/N_1\cap AN_2=H/D\in \mathfrak X.
$$
By the inductive assumption, $H\in \mathfrak X$. It follows that $A$ 
is $\mathbb P$-subnormal in $H$, $H$ is $\mathbb P$-subnormal in $G$, and 
hence $A$ is $\mathbb P$-subnormal in $G$. This contradicts the fact that 
$G$ has a primary cyclic subgroup which is not 
$\mathbb P$-subnormal in $G$. Thus $\mathfrak X$~is formation.

We prove that $\mathfrak X$ is a saturated formation by induction on $|G|$ . 
Suppose that $\Phi (G)\ne 1$ and $G/\Phi (G)\in \mathfrak X$. Since by Lemma 13, 
the quotient group $G/\Phi (G)$ possesses an ordered Sylow tower of 
supersolvable type, it follows that $G$ possesses an ordered Sylow tower of 
supersolvable type.

Let $N$ be a minimal normal subgroup of $G$. It is clear that  
$$
\Phi (G)N/N\subseteq \Phi (G/N), \
G/\Phi (G)N\simeq \big {(}G/\Phi (G) \big {)} \big {/} 
\big {(}\Phi (G)N/\Phi (G) \big {)}\in \mathfrak X,
$$
$$
(G/N)/(\Phi (G/N))\simeq \big {(}(G/N)/(\Phi (G)N/N)\big {)}\big {/}
\big {(}(\Phi (G/N)/(\Phi (G)N/N)\big {)}\in \mathfrak X,
$$
because 
$$
(G/N)/(\Phi (G)N/N)\simeq G/\Phi (G)N\in \mathfrak X.
$$
By the inductive hypothesis, we have $G/N\in \mathfrak X$. Since 
$\mathfrak X$ is a formation, this implies that $N$~is  
a unique minimal normal subgroup of $G$, $N\subseteq \Phi (G)$, 
$N$~is a $p$-subgroup for the largest  
$p\in\pi (G)$, and $O_{p^{\prime}}(G)=1$. Let $P$~be a Sylow 
$p$-subgroup of $G$, $P$ is normal in $G$.

Suppose that $G$ has a primary cyclic subgroup $A$ which is not 
$\mathbb P$-subnormal in $G$. Since $G/N\in \mathfrak X$, it folows that 
the quotient group $AN/N$ is $\mathbb P$\nobreakdash-\hspace{0pt}subnormal in $G/N$,
and by Lemma 3~(2), $AN$ is $\mathbb P$-subnormal in $G$. 
By Lemma 3~(3), we see that $A$ is not $\mathbb P$-subnormal in $AN$, and 
Lemma 5 implies that the orders of $A$ and $N$ are coprime. Therefore, $AP/N$~is
a biprimary subgroup in which the Sylow subgroups $AN/N$ and $P/N$ 
are both $\mathbb P$-subnormal. Theorem 2.13~(2) \cite{VVTSMJ} implies 
that $AP/N$ is supersolvable. By Lemma 15,  
$\Phi (P)=P\cap \Phi (G)$, thus $N\subseteq \Phi (P)\subseteq \Phi (AP)$,
and by Theorem VI.8.6 \cite{Hup}, we deduce that $AN$ is supersolvable. 
Lemma 2~(2) implies that $A$ is $\mathbb P$-subnormal in $AN$, which is a 
contradiction.  

3. Let $G\in \mathfrak X $ and $B$ is a biprimary group with cyclic Sylow
subgroup $R$. By Lemma 16\,(1), $G$ possesses an ordered Sylow tower of 
supersolvable type. If $R$ is normal in $B$, then $B/R$ is primary, it follows
that $B$ is supersolvable. 
If $R$ is not normal in $B$, then $B=PR$, where $P$~is a normal Sylow subgroup
of $B$. By hypothesis, we conclude that $R$ is $\mathbb P$-subnormal in $G$, 
and by Lemma 4\,(1), $R$ is $\mathbb P$-subnormal in $B$. 
Hence, by Lemma 10\,(2), $B$ is supersolvable.

Conversely, suppose that $G$ possesses an ordered Sylow tower of supersolvable type
and every biprimary subgroup of $G$ with cyclic Sylow subgroup is supersolvable.
Assume that $G\not\in\mathfrak X$. Let us choose among
all such groups a group $G$ with the smallest possible order.
Then $G$ contains a cyclic non-$\mathbb P$-subnormal $r$-subgroup $R$.  
Since $G$ possesses an ordered Sylow tower of supersolvable type, we deduce that 
a Sylow $p$-subgroup $P$ for the largest prime $p\in \pi (G)$ is normal.
If $p=r$, then $R\subseteq P$, it follows that $R$ is $\mathbb P$-subnormal in $G$, 
which is a contradiction. Thus $p\ne r$ and $PR$~is biprimary with cyclic Sylow
subgroup $R$. By hypothesis, $PR$ is supersolvable, and by Lemma 2\,(2), $R$ 
is $\mathbb P$-subnormal in $PR$. 
The quotient group $G/P$ possesses an ordered Sylow tower of supersolvable type
and every its biprimary subgroup with cyclic Sylow subgroup is 
supersolvable. Thus $G/P \in \mathfrak X$. 
Since $PR/P$~is a cyclic $r$-subgroup, we see that $PR/P$ 
is $\mathbb P$-subnormal in $G/P$. It follows by Lemma 3\,(2) that $PR$
is $\mathbb P$-subnormal in $G$. Now by Lemma 3\,(3), we obtain that $R$ 
is $\mathbb P$-subnormal in $G$. This is a contradiction. The assertion is proved.

4. Let $G\in \mathcal M(\mathfrak X)$, and let $q$ be the smallest prime 
divisor of $|G|$. Consider an arbitrary proper subgroup $H$ of $G$.
Since $H\in \mathfrak X$, so by Lemma 13~(1), the subgroup $H$ has an 
ordered Sylow tower of supersolvable type, 
in particular, $H$ is $q$-nilpotent. By Theorem IV.5.4 \cite{Hup}, the group
$G$ is either $q$-nilpotent or a $q$-closed Schmidt group. If $G$
is a $q$-closed  Schmidt group, then $G$ is a biprimary  minimal 
non-supersolvavle group whose non-normal Sylow subgroup is cyclic.
In this case, the statement is true.

Suppose that $G$~is a $q$-nilpotent group. Then 
$G=[G_{q^{\prime}}]G_q$. Since $G_{q^{\prime}}\in \mathfrak X$, 
it follows by Lemma 13~(1), that  $G_{q^{\prime}}$ possesses an 
ordered Sylow tower of supersolvable type, and thus $G$ possesses 
an ordered Sylow tower of supersolvable type. 
Let $N$~be a minimal normal subgroup of $G$. 

First, assume that $\Phi (G)=1$. In this case, $G=[N]M$, where $M$~is some
maximal subgroup of $G$. Since $G\not \in \mathfrak X$, then $G$ contains
a primary cyclic non-$\mathbb P$-subnormal subgroup. Let $A$ be a subgroup of 
least order among these subgroups. Since
$$
AN/N\subseteq  G/N\simeq M\in \mathfrak X,  \ AN/N\simeq A/A\cap N,
$$
it follows that $AN/N$ is $\mathbb P$-subnormal in $G/N$, and by Lemma 3~(2),
the subgroup $AN$ is $\mathbb P$-subnormal in $G$.  If $AN\ne G$, then 
$AN\in \mathfrak X$,  it follows that $A$ is 
$\mathbb P$-subnormal in $AN$, and by Lemma 3~(3), the subgroup $A$ is 
$\mathbb P$\nobreakdash-\hspace{0pt}subnormal in $G$, which is a contradiction.
Therefore $AN=G$, in particular, $G$ is biprimary. Let $H$~be an arbitrary 
maximal subgroup of $G$. Then either $A^x\subseteq H$, $x\in G$ or
$N\subseteq H$. If $A^x\subseteq H$, then $A^x=H$ because $N$~is a minimal
normal subgroup of $AN=G$ and $H$ is cyclic. If $N\subseteq H$, then
by the Dedekind identity, $H=(A\cap H)N$. By the choice of $A$ we can 
conclude that $A\cap H$ is $\mathbb P$-subnormal in $G$. Now by Theorem 
2.13~(2) \cite{VVTSMJ}, $H$ is supersolvable. So in the case of $\Phi (G)=1$, 
we proved that $G$ is a biprimary  minimal non-supersolvavle group.

Let $\Phi (G)\ne 1$. According to statement 2 of the theorem, 
$G/\Phi(G) \not \in \mathfrak X$, and so by the inductive hypothesis, 
$G/\Phi(G)$ is a biprimary  minimal non-supersol\-vable group.  
It follows from the structure of such groups that
$G/\Phi (G)$ possesses an ordered Sylow tower. 
Since $\pi (G)=\pi (G/\Phi (G))$, we deduce that $G$~is a 
biprimary group which possesses an ordered Sylow tower: 
$G=[P]Q$, where $P$ and $Q$~are Sylow subgroups of $G$. Since 
$G\not \in \mathfrak X$, then there exists a primary cyclic 
non-$\mathbb P$-subnormal subgroup.  Let $A$ be a subgroup of 
least order among these subgroups. If $PA\ne G$, then $PA\in \mathfrak X$, 
and thus $A$ is $\mathbb P$-subnormal in $PA$. Since $PA$ is $\mathbb P$-subnormal 
in $G$, it follows by Lemma 3~(3) that $A$ is 
$\mathbb P$\nobreakdash-\hspace{0pt}subnormal in $G$, this is a contradiction. 
Therefore, $PA=G$. Let $H$~be an arbitrary maximal subgroup of $G$. Then either 
$P\subseteq H$ or $A^x\subseteq H$, $x\in G$. 
If $P\subseteq H$, then $H=[P](A\cap H)$ by the Dedekind identity. 
By the choice of $A$, we can conclude that $A\cap H$ is $\mathbb P$-subnormal in 
$G$. Now by Theorem 2.13~(2) \cite{VVTSMJ}, the subgroup $H$ is supersolvable.
If $A^x\subseteq H$, then $H=[P\cap H]A^x$. Since $H\in \mathfrak X$, we 
deduce that $A^x$ is $\mathbb P$-subnormal in $H$, and thus $H$ is supersolvable
by Theorem 2.13~(2) \cite{VVTSMJ}. Hence, in the case of $\Phi (G)\ne 1$, 
we proved that $G$ is a biprimary  minimal non-supersolvavle group. 

Therefore, in any case every minimal non-$\mathfrak X$-group is a biprimary
minimal non-supersolvable group. We conclude by Lemma 14, that every non-normal 
Sylow subgroup of $G$ is cyclic. 

The theorem is proved.

\label{kon}

\noindent V.\,N. KNIAHINA

\noindent  Gomel Engineering Institute,\\ 
Rechitskoe Shosse 35a,\\
Gomel 246035, BELARUS

\noindent E-mail: knyagina@inbox.ru

\bigskip

\noindent V.\,S. MONAKHOV

\noindent Department of mathematics\\ 
Gomel F. Scorina State University,\\ 
Sovetskaya str. 104,\\
Gomel 246019,  BELARUS 

\noindent E-mail: Victor.Monakhov@gmail.com

\end{document}